# Rolle Theorem and Bolzano-Cauchy Theorem from the end of the 17th century to K. Weierstrass' epoch

G.I. Sinkevich,

Saint Petersburg State University of Architecture and Civil Engineering

*Under consideration is the history of a famous Rolle's theorem as follows: "If a function is continuous at [a, b], differentiable in (a, b), and $f(a)=f(b)$, then at least one such point c will be found in (a, b) that $f'(c)=0$", and the history of the related theorem on the root interval as follows: "If a function is continuous at [a, b] and has different signs at the ends of the interval, then at least one such point c will be found in (a, b) that $f(c)=0$". In the 20th century, this theorem became known as Bolzano–Cauchy Theorem.*

Speaking of the 19th century reform of analysis, we recollect its key characters, in the first place. They were A. Cauchy and K. Weierstrass. In fact, quite a lot of scientists form part of its real history. Some of them lived earlier (and therefore, in many ways contributed to the development of analysis), others were their contemporaries without whom this history would lack not only certain important colors, but the thoroughness a historical evidence needs as well.

Mathematician and theologist Bernard Bolzano was the ingenious messenger of the ideas of this reform. Many fundamental ideas the implementation whereof is normally associated with the names of the abovementioned Cauchy and Weierstrass, as well as G. Cantor and R. Dedekind, belong to him. These are the notion of the least upper bound (1817), realization of the need in and attempt to develop the theory of a real number (1830s), set-theoretic understanding [1].

In 1817, Bernard Bolzano wrote a work entitled "Purely analytic proof of the theorem that between any two values which give results of opposite sign, there lies at least one real root of the equation" [1, 43].

Bolzano attributed the importance of the key property of a continuous function to this theorem and considered its genesis. Let us follow his lead.

Before the 17th century, they used geometric images of curve crossings to find roots of algebraic equations, while the interval in which the roots lied was determined on the basis of ratio analysis. For example, Newton wrote in his *Universal Arithmetics*: "Should you want to find a limit which cannot be exceeded by any root, find a sum of squared roots and take the square root of this sum. This square root will be larger that the largest equation root" [11, p. 265]. In the middle of the 17th century, they used the method of root localization with the help of an auxiliary equation. As a rule, they would search for positive roots. In order to make an auxiliary equation, exponents of variables were reduced by a unity, and each factor was multiplied by the former exponent. (Johann Hudde, 1658, – J. Hudde, 1628–1704)[2]. Later, this operation was defined as differentiation of a polynomial (Isaak Newton and Gottfried Leibnitz).

## 1690, Rolle and his Method of Cascades

Michel Rolle (M. Rolle, 1652–1719) was born in France in a small town of Ambert, province Auvergne, to a family of a shoemaker. He came to Paris when he was 23 to work there as an enumerator. He achieved such success in self-education that in 1682 managed to solve a complicated problem of Jacques Ozanam (1640-1717): "Find four numbers such that the difference of any two is a square, and the sum of any two of the first three is still a square[3]". Ozanam himself believed that each number consists of at least 50 digits, however, Rolle found such numbers, each containing no more than 7 digits. By this solution he gained a mathematical reputation. He was invited to teach the son of the Minister of War, got an appointment in the War Ministry, a pension from Louis IV, and in 1685, became member of the Royal Academy of Science (as a student of an astronomer), and from 1699, its welfare recipient (as a geometrician, that is to say, mathematician).

---

[1] Bolzano's ideas gained popularity in Germany thanks to Hermann Hankel who in 1870 published the abovementioned Bolzano's work in Tubingen and popularized his other works. Otto Stolz' work devoted to Bolzano [42] is of interest as well. Cauchy knew Bolzano personally and used his ideas in works of his own [27].
[2] The reconstruction of Hudde method was provided by A.P. Yushkevich in his comments to the translated work of L'Hôpital entitled "Analysis of Infinitely Small" [2, p. 400].

Rolle dealt with algebraic issues: Diophantine analysis, solving algebraic equations. He largely popularized R. Recorde's algebraic symbolism and introduced a symbol $\sqrt[n]{x}$.

Rolle is known for his violent criticism of differential calculus and Method of Descartes for the lack of adequate rationale. French mathematicians P.Varignon and J. Saurin refuted most Rolle's arguments; in 1705, the Academy recognized that he was wrong, to which he later agreed. However, this discussion compelled Leibnitz to state differential calculus more strictly.

In 1690, Rolle published *A Treatise on Algebra* [3] devoted to solving of Diophantine and algebraic equations with arbitrary exponents. The statement contained many new ideas which were innovative compared to the same Method of Descartes. Method of Cascades as one of the methods [3, p. 124–152] was based on the idea that roots of the initial equation have been divided by roots of an auxiliary derived equation. Roots of the auxiliary equation can also be separated with the help of another auxiliary equation, etc. Forming cascades, we descend to a linear equation; and once we solve it, we ascend back to the initial one. Rolle published the justification of his method a year later in a small work entitled "Démonstration d'une méthode pour résoudre les égalités de tous les degrés" (Justification of the solving method for equations with whatsoever exponents), Paris, 1691.

---

[3] Trouver quatre nombres tels que la différence de deux quelconques soit un carré, et que la somme de deux quelconques des trois premiers soit encore un carré. In English literature this text was translated without «of any two".

Let us demonstrate how Rolle himself solved a quartic equation from his Treatise [3, p. 124-126]:

$$v^4 - 24v^3 + 198vv - 648v + 473 \infty 0,$$

(in Rolle's record, sign $\infty$ stands for equality, sign $vv$ stands for $v^2$).

In order to find a root interval, he introduced a Grande Hypothèse (great hypothesis), a Petite Hypothèse (small hypothesis), and a Hypothèse extrême (extreme hypothesis). The great hypothesis was defined as the value for which no roots of a polynomial can exceed which was calculated as follows: $\left(\dfrac{a}{c}+1\right)$, where $a$ was the absolute value of the largest negative coefficient. Here, $a = |-648| = 648$; $c$ was a coefficient at the most significant exponent, $c = 1$, therefore, the great hypothesis for our equation equaled 649. Rolle asserted that this was true for all polynomials. The small hypothesis was the number that is less than every root. Whereas only positive roots were considered, it was therefore zero which was taken as the small hypothesis. The extreme or end hypotheses were the intermediate bounds dividing roots of auxiliary equations known as cascades. Rolle called all extreme hypotheses 'Hypothèses moyennes' (intermediate hypotheses).

Thus, $v^4 - 24v^3 + 198v^2 - 648v + 473 = 0$, all positive roots lied on the interval (0, 649). Rolle's task was to split this interval into several smaller intervals, each containing only one root of the initial equation, i.e. select limits of root (root interval).

The first component of the sum had the exponent of the unknown which equaled 4. Rolle multiplied it by 4; the exponent of the second component of the sum was 3, and he multiplied it by 3. The exponent of the third component of the sum was 2, and he multiplied it by 2; the exponent of the fourth component of the sum was 0, and he multiplied it by 0. Further, he had $4v^4 - 72v^3 + 396v^2 - 648v = 0$ and divided all members of the equation by the unknown $v$, whereupon, he obtained $4v^3 - 72v^2 + 396v - 648 = 0$. Again, he multiplied each member of the equation by the exponent: $12v^3 - 144v^2 + 396v = 0$ and divided by the

unknown: $12v^2 -144v+396=0$. And again: $24v^2 -144v=0$, $24v-144=0$, $4v-24=0$.

Let us locate the cascades as follows:

The first cascade: $4v-24=0$.

The second cascade: $6v^2 -72v+198=0$.

The third cascade: $4v^3 -72v^2 +396v-648=0$.

The fourth cascade: $v^4 -24v^3 +198v^2 -648v+473=0$.

Roots of each cascade are divided by roots of the previous one, and all positive roots are lying within (0, 649). Whereas 6 is the root of the first cascade, the roots of the second cascade therefore is lying within (0, 6) and (6, 13), where each interval contains only one root, and number $13=\frac{|-144|}{12}+1$ is the great hypothesis for this equation (all symbols are up-to-date[4]). We will be only finding the leftmost root, the rest of them to be calculated in the same way. The values of polynomial $6v^2 -72v+198$ on the bounds of the interval (0, 6) have different signs. Let us take any average value from the interval, not necessarily from the middle, and check the signs:

$$f_2(5)=-12<0, \; f_2(4)=6>0.$$

Accordingly, the root of the second cascade lies within (4; 5). Proceeding with iterations or using the well-known Viet formula, we will obtain value $6-\sqrt{3}$. The second root is $6+\sqrt{3}$.

Accordingly, the bounds within which roots of the third cascade lie will be as follows: $(0; 6-\sqrt{3}), (6-\sqrt{3}; 6+\sqrt{3}), (6+\sqrt{3}; 163)$; besides, there is only one root lying in each interval. Here, value $163=\frac{|-648|}{4}+1$ is the great hypothesis for the third cascade. We will be only finding the leftmost root of the third cascade. Let us check signs of the third cascade at the ends of this interval. They are different. Let us take any middle point from interval $(0; 6-\sqrt{3})$ and calculate the value of $f_3(v) = 4v^3 - 72v^2 + 396v - 648$.

---
[4] The modulus sign was introduced by K. Weierstrass.

$f_3(0) = -648 < 0, f_3(5) = 32 > 0, f_3(4) = 40 > 0, f_3(3) = 0.$

Accordingly, the left root of the third cascade equals 3, which means that the left (positive) root of the fourth cascade (i.e. root of the initial equation) lies within (0, 3). Let us calculate $f_4(v) = v^4 - 24v^3 + 198v^2 - 648v + 473$ at the ends of the interval and in certain middle points: $f_4(0) = 473 > 0, f_4(3) = -256 < 0, f_4(1) = 0.$

By doing so, we found the left root of equation $v = 1$. All other roots will be calculated in the same way. If we are doing an approximate calculation, the procedure enables us to find the root with an accuracy to any decimal sign. For long intervals, Rolle applies an auxiliary equation replacing $v = \left(\frac{a}{c} + 1\right) - x$, where $\left(\frac{a}{c} + 1\right)$ is the great hypothesis.

In addition to this method for solving algebraic equations, Rolle offers four more methods as well methods for solving indefinite equations, and a method of finding the polynomials common divisor.

As we can see, using coefficients, Rolle selected limits between which the roots were lying. In the method of cascades, although without differential calculus terminology, he used the principle of partitioning of root of the polynomial by roots of its derivative, and the existence of roots was checked by the difference of polynomial's signs at the ends of the interval. In 1691, in his work devoted to the justification of the method of cascades, Rolle demonstrated that the values of the derivative (i.e. the derived polynomial) for two adjacent (single) roots of the integral polynomial have different signs [4, p. 47].

The notion of a function was only evolving in the 17$^{th}$ century; there was no notion of a function graph or geometrical locus at that time. Therefore, the concept of the root as a point where the function graph crosses the axis did not exist as yet.

This image appeared in Michel Rolle's works. He verified the existence of a root in the interval determining signs of the polynomial in the left part of the equation at the ends of the interval. If the signs were different, the root was lying within the interval. Rolle narrowed the interval checking the sing of the polynomial in any inner point of the interval. Thus, Michel Rolle became the father of two analysis theorems: the theorem on the root interval currently known as Bolzano-Cauchy Theorem and Rolle Theorem as such stating that roots of a continuous function are divided by roots of its derivative.

In Russian, this method was stated in a work of S.A. Yanovskaya [5]; there is a good reconstruction of the method of cascades in English [6]. The primary source of Rolle's biography is *Eloge de M. Rolle* written by his contemporary, Secretary of Paris Academy of Sciences B. Fontenelle [7]. In more detail, the biography, Rolle's method and its history are laid down in an article by G. Sinkevich [8].

**1707, Rolle and Newton**

Before Rolle, approximate solution of algebraic equations was achieved by graphic methods, i.e. by curve crossing and with the help of simple iterations. Rolle was probably the first to form the notion of the root interval by comparing signs of the respective polynomial.

The very geometric image of the problem did not correspond to the search of the curve crossing with the axis. Instead, it was the search of the points of intersection of two curves. Therefore, an image of a graph with ordinates of different sign at the ends of the interval could not appear in terms of algebra. I. Newton (1642–1727), for example, imagined variables to be time variant, not changing against each other

[9]. The concept of a line as a geometrical locus of the equation[5] first appeared in the work of L'Hôpital devoted to conic sections; thereafter, it was developed by Euler, and the general approach formed only in the 19$^{th}$ century.

Newton described his tangent method in his works entitled "Analysis of equations with an infinite number of members[6]" and "Method of fluxions and infinite series[7]". This method was also stated in the book of 1685 by J. Wallis *A Treatise of Algebra both Historical and Practical*. In 1690, J. Raphson's (1647/48–1715) treatise was published in England. It was entitled "Analysis aequationum universalis" [10] and contained an improved statement of Newton-Raphson method or tangents method[8]. In 1707, Newton's *Arithmetica Universalis* (Universal Arithmetics) was published. It contained equations numerical solution methods [11].

Using the tangents method, Newton did not check signs of the function at the ends of the interval until this Rolle's treatise appeared, which appears from Newton's works published before 1690, e.g. in his *Method of Fluxions* of 1671 [12]. To determine the starting point for the root calculation procedure, Newton used the method of false position and the so-called 'Newton parallelogram' or 'Newton polygon'.

Paris Academy and London Royal Society exchanged academic literature. No doubt, Newton received Rolle's Treatise. Moreover, he included the statement of his method in his 1707 publication of Universal Arithmetics [11, p. 267–270] without mentioning Rolle's authorship though. But it was only after 1707, when Rolle's work appeared, that Newton first started checking signs of polynomials at the ends of the interval in his *Universal Arithmetics*. The method of narrowing the interval containing the root by checking the polynomial sign in a certain inner (not necessarily middle) point first occurred in Rolle's works. Bolzano formalized it 117 years later as half-interval method. Please note that Newton believed all functions under consideration to be defined as continuous, while Rolle considered only polynomials which are continuous functions.

---

[5] This concept should be distinguished from that of a geometrical locus of points which possess a given property (e.g. circumference) which appeared in the ancient world.
[6] "De analysi per aequationes numero terminorum infinitas", manuscript of lectures Newton read at the university, was written in Latin in 1669 and published in 1711.
[7] "De methodis fluxionum et serierum infinitarum", 1671, translated into English and published as "Method of Fluxions" in 1736.
[8] While Newton considered the sequence of approximating polynomials, Raphson already considered successive iterations of the variable.

The Newton method and use of Maclaurin series expansion were more popular among continental mathematicians. In 1740, T. Simpson provided a summary of the Newton method in his work entitled "Essay on several subjects in speculative and mixed mathematics[9]" [13].

### 1696, L'Hôpital

In 1696, the first textbook on analysis was published in Paris. It was *Analyse des Infiniment Petits pour l'Intelligence des Lignes Courbes* [14] by marquis G.F. de L'Hôpital (l'Hospital, 1661–1704) stating lectures of I. Bernoulli (1667–1748). It was for the first time that differential and integral calculus was stated there and the notions of abscisse, ordinate, coordinates, and geometrical locus were used. However, L'Hôpital achieved the simplicity and intelligibility of presentation by neglecting argumentation: "I am sure that in mathematics only conclusions are of use and that books describing but details or particular suggestions only make those who write them and those who read them waste time." [2, p. 57] But he provided the geometrical meaning of a derivative, relation between the function increase and decrease and the sign of the first derivative, the necessary condition of the extremum. One would also find the following reasoning there: "No continuously increasing or decreasing positive value[10] can turn into negative without passing through infinity or zero, more specifically: through zero when it is first decreasing and through infinity when it is first increasing. This implies that the differential of the largest or the least value must equal zero or infinity. It is easy to understand that a continuously decreasing positive value cannot turn into negative

---

[9] Simpson already used derivatives for iterations.
[10] L'Hôpital means subtangent.

without passing through zero; however, it is not so obvious that in the case of an increase it has to pass through the infinity." [2, p. 130–132]

This book opens the initial period of development of analysis where all functions were continuous as they were algebraic whole numbers and all analytical statements were based on geometrical ideas. The rules of differential calculus of the 17th-18th century were defined for algebraic functions only. Formulas of derived transcendental functions will appear later in Euler's and Cauchy's works.

**1708, Rolle's Method in treatise by Reynaud**
C.-R. Reynaud (1656–1728), French preacher and philosophy professor, was familiar with works of Hudde, Descartes, Rolle, Newton, Leibnitz, Bernoulli, and L'Hôpital. Being aware of the reproaches in relation to the insufficiency of argumentation and lack of methodical statement of modern mathematics, he made it his crusade to deliver a complete course of analysis, algebra and geometry with reference to one another and demonstrations. In 1708, his book entitled "Demonstrable analysis" [15] in two volumes was published in Paris stating results of the above mathematicians. I would note that before that time, only geometrical statements were demonstrated in mathematics, while in algebra and evolving analysis they only showed examples.

The first volume of *Demonstrable Analysis* was devoted to algebraic issues and the second one, to differential and integral calculus where the author tried to demonstrate most of the statements, provided quite a lot of examples, not only mathematical ones but from mechanics and astronomy as well. Notes to the second edition of 1736 were written by Varignon.

Professional preacher, Reynaud had good command of presentation and was good at selecting terms not only in Latin but in French as well. His manner of presentation compared favorably with the complicated language of Michel Rolle not only in terms of speech melody in general but in terms of consistency of argumentation and appropriateness of definitions as well. One can feel an experienced teacher in his techniques. First, Reynaud considered linear equations forming equations of the highest degrees by multiplying binomials, proceeding up to equations of the sixth degree. He showed solutions to not only numerical but to algebraic letter equations as well, both in radicals and approximately; he provided the fundamental theorem of algebra – the theorem on the value of the residue of division of a polynomial by a binomial[11] [15, v. I, p. 270–271]. His demonstrations included verbal proof accompanied by demonstration of particular cases and examples. In algebraic equations, Reynaud distinguished between cases with single, multiple, positive, negative, whole, fractional, non-measurable, and imaginary roots.

Reynaud devoted a large section [15, v. I, p. 269–375] to the method of Rolle: "The sixth book explains and proves the method of finding the values which constitute limits of the unknown in the exponential equation (Mr. Rolle was the author of this method) and provides some solution methods using these limits; the roots here can be found with any degree of accuracy as may be wished" [15, v. I, p. XII].

In development of Rolle's ideas, Reynaud introduced terminology of his own. He called each auxiliary equation the roots whereof are limits of roots of the previous one 'l'équation des limites' (equation of limits) and boundaries of the interval which contains a root, 'limits of roots'. Reynaud defined the root interval which contains a real root based on the difference of signs in the left part of the equation on the limits of the roots; he described Rolle's step-up and step-down procedures.

The second volume of *Demonstrable Analysis* is devoted to differential and integral calculus. It contains a statement that a tangent to a curve (for conic sections) in a certain point is parallel to diameter [15, v. 2, p. 176). Another statement it contains reads as follows: if a sequence of values of a variable (e.g. subtangent[12]) is first positive and thereafter, becomes negative, this means that it passes a certain point where its value equals zero or infinity [15, v. 2, p. 177].

---

[11] This theorem bears the name of E. Bezout (1730-1783). One can come across this very theorem in J. Raphson's (1648-1715) works.

### 1727–1729, Rolle Theorem in Campbell's and MacLaurin

Reynaud's *Demonstrable Analysis* was well known in England. G. Campbell[13] quoted it in 1727 in his work entitled *A Method for Determining the Number of Impossible[14] Roots in Adfected Aequations* [16]. Campbell translated French mathematical works into English and solved algebraic equations himself. In the above work, he paraphrased Rolle's procedure, engaged Fermat's rule for determining maxima and minima, and considered the case with the final quadratic equation with a negative discriminant. Then, judging from the change of signs of coefficients of the initial equation, it is possible to tell the number of imaginary roots or rather their least number. In his letter [17] to the same journal, C. MacLaurin (1698–1746) argued against rigidity of the demonstration with him. MacLaurin formulated the theorem as follows: "Roots of equation $x^n - Ax^{n-1} + Bx^{n-2} \& c. = 0$

are limits of roots of equation
$nx^{n-1} - (n-1)Ax^{n-2} + (n-2)Bx^{n-3} \& c. = 0$"

[17, p. 88]. The change of signs of the coefficients ensures *n* positive roots. This statement had already had the status of a theorem. MacLaurin also considered a more general type of an auxiliary equation used by Hudde in 1658.

### 1746, A.C. Clairaut

In 1746, a book of A.C. Clairaut (1713–1765) *Fundamentals of Algebra* [18] was published to tell about methods to solve algebraic equations. Much attention was given to the method of Newton and MacLaurin there, however, Rolle and his method were not mentioned.

---

[12] A subtangent is a projection of an interval of a tangent between the points of intersection with axis *OX* and tangency point on axis *OX*.
[13] This Scottish mathematician is only known to argue against C. MacLaurin and die in 1766.
[14] Imaginary ones.

## 1755, Rolle Theorem in Euler

In 1755, St. Petersburg Academy of Sciences published L. Euler's (1707–1783) work entitled "Institutiones calculi differentialis" (Foundations of Differential Calculus). The algebra and analysis convergence trend we saw in Reynaud's treatise and discussions of the equation of string by D'Alembert and Euler led to an extension of the notion of a function. Euler was proud that he did not have to turn to applied interpretation when stating analysis. He wrote in Chapter IX: "The notion of the equation can be traced to the notion of function" [19, p. 367]. Euler considered a polynomial there as a trivially continuous function which satisfied his concept of a continuous function as a function assigned by an integrated analytic expression. Euler repeated the above theorem of MacLaurin on roots of equation

$$x^n - Ax^{n-1} + Bx^{n-2} - Cx^{n-3} + dx^{n-4} - etc... = 0,$$

divided by roots of the auxiliary equation, i.e. by extreme values. Please note that both Maclaurin and Euler considered an equation which trivially had $n$ real roots. The problem of determination of the number of imaginary roots set up by Campbell was considered by Euler in Chapter XIII. Euler summarized his reasoning as follows:

"However, it appears from the above that, although not all roots of the proposed equation may be real, nevertheless, there is always a maximum and a minimum between any two roots. As to the converse proposition, it is wrong in general, i.e. there may be no real root between any two maxima or minima. However, this conclusion may be made subject to an added condition that either value of $z$ will be positive and the other one, negative […]. There is one value between two real roots of the equation at which the functions becomes a maximum or minimum." [19, p. 435–436] Euler's reasoning was based on the concept of continuous movement; he extended all properties of algebraic expressions to functions.

### 1797, Theorem on the Root Interval in S.F. Lacroix' Élémens d'algèbre

In 1797, the first edition of S.F. Lacroix' (1765–1843) *Élémens d'algèbre* was published in Paris. Lacroix was the author of courses of higher mathematics which were repeatedly reissued and well known in Russia of the 19$^{th}$ century.

In the *Élémens d'algèbre*, he provided the following theorem on the root interval:

"If there are two values which, being inserted into the equation instead of an unknown, will produce two results opposite in sign, we can conclude that the roots of this equation are between these two values and they are real." [20, p. 298].

In 1811, an authorized German translation of Lacroix' *Élémens d'algèbre* by M. Metternich (1747–1825), professor of mathematics and physics of Mainz University [21], was published in Mainz. This book also contains the statement of the theorem on the root interval. This book was repeatedly reissued in German and widely used by German mathematicians.

### 1768, Kaestner about selection of the root interval

Abraham Gottheld Kästner (1719–1800), professor of mathematics and physics in Göttingen, was esteemed as a prominent teaching methods specialist on various issues relating to analysis. Please note that he considered irrational numbers as limits of rational sequences before Cauchy. Kaestner was in correspondence with Euler.

In 1768-69, Kaestner wrote a course of *Fundamentals of Mathematics* (Der mathematischen Angfangsgründe) in four volumes [22]. It was an accomplished course in terms of methodology; it provided a good historical survey; and was repeatedly reissued. One can feel Euler's influence in the course. The course of Kaestner was published in Russian in 1792–1803. There was no Rolle's method in the course of Kaestner, but there is the theorem on root interval of a polynomial with a demonstration using geometric analogy. In the German publication of 1794, it reads as follows: "Theorem. If $y$ is positive for $x = a$ and negative for $x = c$, then at least one value of $x = b$ will be found between $a$ and $c$ for which $y = 0$" [22, p. 198]. Kaestner has great influence on German mathematical education, Karl Weierstrass addressed his works.

**1798, Lagrange about the Method of Rolle**

In 1798, J.L. Lagrange (1736–1813) suggested his root isolating method based on the method of Rolle [23]. Lagrange asserted that roots of the initial equation were divided by roots of the derived equation and characterized by insertion of roots of the derived equation in the initial equation followed by determination of its sign: "Thus, these rules enable us to determine not only the number of real roots of the equation, but the boundaries within which they lie as well; and if you wish to constrain roots between a value that is larger than $\alpha_1$ and less than $v_1$, an additional search needs to be carried out in accordance with the method stated in Chapter IV (No. 12) regarding the boundaries of positive roots of this equation.

Please note that already Rolle knew the rules enabling us to find these limits and stated by us according to Newton and MacLaurin, which appears from Chapters V and VI of this *Algebra*." [23, p. 199].

**1817, Bolzano and the theorem on root interval**

Bernard Bolzano (1781–1848), Czech mathematician and philosopher, contributed quite a lot in the development of the notion of continuous and infinite. In his manuscript of 1817 entitled "Purely analytic proof of the theorem that between any two values which give results of opposite sign, there lies at least one real root of the equation" [1], he criticized demonstrations of Kaestner, Clairaut, Lacroix, Metternich, Rösling, Klügel, and Lagrange for the involvement of geometrical and physical images (time and movement, transfer) and the lack of analyticity in their reasoning, i.e. lack of understanding of the continuity as a mathematical notion[15]. Bolzano wrote: "As a matter of fact, if we take into account that a proof in science must not at all be just *words* but *argumentation*, i.e. be the exposition of objective cause for the true being proved, then it goes without saying, that if an affirmation is correct only for the values in the space, it may not be correct for *all* variables, whether or not they are in the *space*. The most common kind of proof depends on a truth borrowed from geometry, namely, that every continuous line of simple curvature of which the ordinates are first positive and then negative (or conversely) must necessarily intersect the *x*-axis somewhere at a point that lies in between those ordinates. There is certainly no question concerning the correctness, nor indeed the obviousness, of this geometrical proposition. But it is clear that it is an intolerable offense against correct method to derive truths of pure (or general) mathematics (i.e., arithmetic, algebra, analysis) from considerations which belong to a merely applied (or special) part, namely, geometry. No one will deny that the concepts of *time* and *motion* are just as foreign to general mathematics as the concept of *space*. We strictly require only this: that examples never be put forward instead of *proofs* and that the essence of a deduction never be based on the merely metaphorical use of phrases or on their related ideas, so that the deduction itself would become void as soon as these were changed» [1, italics by Bolzano, p. 172–174; 43].

Bolzano stated the law of continuity[16] as follows: "According to a correct definition, the expression that a function $f(x)$ varies according to the law of continuity for all values of $x$ inside or outside certain limits means just that: if $x$ is some such value, the difference $f(x + w) - f(x)$ can be made smaller than any given quantity provided $w$ can be taken as small as we please» [1, c. 174–175]. The *true* statement is that the continuous function never reaches its top value without first passing through all downs, that is to say that $f(x+n\Delta x)$ can take on any value lying between $f(x)$ and $f(x+\Delta x)$ if $n$ is defined arbitrarily between 0 and 1. However, this statement may not be deemed to be an *explanation* of the notion of continuity, it constitutes a *theorem* of continuity" [1, italics by Bolzano, p. 175]. "Therefore, this theorem can be stated as follows: "If a variable which depends on any other variable $x$ turns to be positive for $x=\alpha$ and negative for $x=\beta$, there is always a value $x$ lying between α and ß for which it becomes zero or a value for which it becomes continuous" [1, p.176]. This is a theorem, said Bolzano, and it must be proved. Bolzano also noted that such point $x$ need not be the only one. He believed it to underlie the algebraic theorem on factorization of a polynomial and Lagrange's theorem on positivity of a definite integral for positive function which equals zero only in the end point of the interval.

---

[15] To tell the truth, as we could satisfy ourselves, neither Clairaut in the above work of 1746 nor Dr. Ch. L. Rösling (1774–1836) from Erlangen University in his book of 1805 entitled "Fundamentals of the theory of forms, differentials, derivatives and integrals of functions" [24] turned to the method of Rolle and his theorem on the root interval. References to them provided by Bolzano are at variance with the topic concerned [1, p.171].
[16] The wording of this law belongs to Leibnitz.

Bolzano suggested his own more stringent plan of proof of this theorem based on another more general one: "If two functions $x$, $f(x)$ and $\varphi(x)$ either for *all* values of $x$ or for all values lying between α and ß change *in accordance with the law of continuity*; if further $f(\alpha) < \varphi(\alpha)$ and $f(\beta) > \varphi(\beta)$, each time there is value $x$ lying between α and ß for which $f(x) = \varphi(x)$" [1, p. 170 – 204, p. 198]. Bolzano proved this theorem assuming that there exists an upper bound of an area where an abstract property of the function[17] is met and using the interval bisecting method in the auxiliary theorem. He demonstrated that the existence of the least upper bound does not cause any inconsistency, as it became possible to provide a stronger demonstration of the existence of a bound only after Weierstrass' works in 1860-th and after 1872 when the theory of a real number was developed by C. Meray, K. Weierstrass, E. Heine, R. Dedekind, and G. Cantor. Bolzano tried to develop the theory of real numbers through section later, in 1830s [25].

---

[17] e.g. negativeness.

Thereafter, Bolzano proved the theorem on the root interval. Here, Bolzano stated another theorem: "Passing from one value to another, at least once, the function takes the value of each intermediate value[18]." Bolzano emphasized that the above property was a result of continuity, however, it could not be taken as basis for defining the continuity.

I would also note that this work contains the sequence convergence criterion [1, 188–189] stated 4 years later by A. Cauchy and named after him.

### 1821, Cauchy, the Theorem on the Root Interval in the "Cours d'analyse de l'École royale polytechnique"

In 1821, A.L. Cauchy (1789–1857) published *The Course of Analysis* [26] lectured in École Polytechnique. The first part of the course was entitled "Analyse algébrique" and its second part, *Le Calcul infinitésimal,* was published in 1823.

Cauchy was brilliant in improving ideas voiced by his colleagues who sometimes failed to provide strong or skilful reasoning. As an example, I would mention ideas of Ampère, Abel, Grassmann, Bolzano, and Galois [27].

The notion of a continuous function introduced in *Algebraic Analysis* exactly repeated Bolzano's definition [26, p. 43].

Cauchy did not mention Rolle, although he addressed the approximate solution of algebraic equations. He provided the following theorem in the chapter devoted to solution of equations [26, p. 378]: "Let $f(x)$ be a real function of the variable *x*, which remains continuous with respect to this variable between the limits $x=x_0$, $x=X$. If the two quantities $f(x_0)$ and $f(X)$ have opposite signs, we can satisfy the equation (1) $f(x)=0$ with one or several real values of *x* contained between $x_0$ and $X$ " [44, p. 330]. Cauchy proved this theorem using interval bisection, but, unlike Bolzano, he did not use the notion of the upper bound. Instead, he was based on convergent sequences.

---

[18] i.e. continuum.

What is very important for analysis, Cauchy stated the theorem on the intermediated value as a property of a continuous function: "*Theorem on a continuous function.* If the function $f(x)$ is a continuous function with respect to the variable $x$ between the limits $x=x_0$, $x=X$, and if $b$ denotes a quantity between $f(x_0)$ and $f(X)$, we may always satisfy the equation $f(x)=b$ by one or more real values of $x$ contained between $x_0$ and $X$" [44, p. 32].

Cauchy provided a couple of methods to solve algebraic equations including the Descartes' method; compared the methods of Newton and Lagrange by the example of solving the same cubic equation; however, he told nothing about dividing roots of the equation by roots of the derivative.

### 1834, Rolle Theorem in Drobisch

In 1834, Professor of Leipzig University M. W. Drobisch (1802–1896) published *The Lectures on Equations of the Highest Orders* where he described Rolle's method of cascades in §107 [28, p. 161]. He complained of Rolle's complicated language, however, called his method worthy of respect, reasoned and stated it as follows: "These theorems (rules) were obtained incomplete in the proposed Rolle's draft. Moreover, it was based on the method of cascades which was extremely hard to understand. There was an essential kernel there that to solve the initial equation, one after another, they formed auxiliary equations, which is much like constructing a house using this method. Thus, we successively obtain roots of low-order equations which provide us with reliable limits of the roots of highest-order equations which we calculate approximately and will describe later up to the roots of the initial equation. This method is based on the assumption that the initial equation has roots in general. This brings about its limitation and impractical awkwardness instead of search for a straighter way." [28, p. 186–188].

Drobisch quoted the theorem on the root interval from Cauchy's course and his mean value theorem in Chapter Seven of *The Alternative Method of Identifying Real and Imaginary Roots* [28, p. 161–176]. He stated the theorem as follows: "Two neighboring real roots are divided by the root of the derivative of an equation, the roots of which derivative are, in turn, divided by roots of the next derivative." [28, p. 176].

"Theorem 1. There is at least one real root of a derivative lying between the two neighboring real roots of the initial equation; however, there may also be 3, 5, and any other odd number of roots between them. Theorem 2. There lies no more than one root of the initial equation between the two neighboring real roots of the derivative equation, however, it may also happen that there are no roots at all between them. Theorem 3. At least one real root of the initial equation may be larger than the largest real root of the derivative equation; no more than one real root of the initial equation may be less than the least root of the derivative equation; however, it may happen that there is no real root of the initial equation which is larger than the largest root and less than the least root of the derivative equation. In this last conclusion, we joined two halves of the $2^{nd}$ conclusion, i.e. the largest root of the initial equation may lie between the first and the second, or between the second and the third, the third and the fourth, etc. root of the derivative equation[19]." [28, p. 178–179].

The algebraic aspect of the theorem of Rolle in solution of equations attracted attention of Italian mathematician G. Bellavitis (1803–1880) who described the method of Rolle in his book entitled "A simple way to find real roots of algebraic equations and a new method of determination of imaginary roots" [29].

### 1861, Rolle theorem in Weierstrass

The vision of continuous functions drastically changed in mid $19^{th}$ century when new mathematical objects appeared, when it was necessary to classify points of discontinuity and assess the scope of this notion and possibility to neglect them when expanding a function in Fourier series. The definition of a continuous function in the language of "$\varepsilon - \delta$" was introduced by Karl Weierstrass (1815–1897) in 1861 [27]; E. Heine (1821–1881), R. Dedekind (1831–1916) and G. Cantor (1845–1918) continued developing the concept of continuity in their works [30, 31, 32] in 1870s.

---

[19] Stated in accordance with Lagrange's Résolution de l'équat. Numér. Not. VIII, who might have probably been the first to make any of the above assertions based on the reputable rule of Rolle. – *Note of W. Drobisch.*

In the summer term (May-June) of 1861, Karl Weierstrass read a course of lectures in differential calculus in the Royal Institute of Trade in Berlin. Herman Schwartz kept his notes of these lectures [33] which were thereafter published by Pierre Dugac [34].

Weierstrass defined the continuous functions and described their properties. "If $f(x)$ is function $x$ and $x$ is a defined value, then, as $x$ passes to $x+h$, the function will change and will be $f(x+h)$; the difference $f(x+h) - f(x)$ is called the change which occurs in the function because the argument passes from $x$ to $x + h$. If such bound $\delta$ can be determined for $h$ that for *all* values of $h$ (the absolute value whereof is still less than $\delta$), $f(x+h) - f(x)$ becomes less than any arbitrary small value $\varepsilon$, then they say that the infinitely small changes in the function correspond to the infinitely small changes in the argument. Because they say that any value can become infinitely small, if its absolute value can become less than any arbitrary small value. If any function is such that infinitely small changes in the function correspond to the infinitely small changes in the argument, they say that this is a *continuous function* of the argument or that it continuously changes together with its argument." [35, p. 189].

There is a theorem as follows under the heading "Study of changes in functions": "*If $f(x_1) = f(x_2)$ for two defined values of $x_1$ and $x_2$ of the argument, then there must be at least one value of $x_0$ between $x_1$ and $x_2$ for which the first derivative $f'(x)$ equals zero*". As A.P. Yushkevich believes, "this is the first or one of the first wordings of the so-called theorem of Rolle" [35, p. 193].

Later, in 1886, analyzing the expansion of the notion of a function, Weierstrass wrote that first, only functions represented by rational number expression, e.g. those next to rational coefficient, were considered. "They changed in accordance with the law of continuity, and this was all we knew about the function. But the discovery of Fourier's series has shown that this is not true; there are continuous functions which may not be obtained by representation as before. There can always be found a mathematical expression for a strictly determined continuous function. Therefore, properties of any function can be derived from the basic notions of continuity, as it is important for any research to derive further notions from the basic ones." [36, p. 21].

The theorem on the root interval, mean value theorem, and theorem on the root of a derivative obtained the status of properties of continuous functions. There were eleven properties like that in Dini's works, while there were about four such properties in the time of Cauchy.

**1878, Rolle Theorem in Dini**

In 1878, a *Course of Lectures on the Theory of Functions of a Real Variable* by U. Dini [37] (1845–1918), professor of Pisa University, was published. It was for the first time that the definition of a continuous function was introduced in this course through unilateral limits. The theorem of Rolle was anonymously stated as follows: "If function $f(x)$ in interval $(\alpha,\beta)$ is finite and continuous in all points but ends of the interval and has a finite and specified or an infinite and specified derivative, and in addition, has the same values in extreme points $a$ and $b$, then there is at least one point $x'$ in interval $(a,b)$ for which point $f'(x')=0$" [37, p. 76–77].

**1879, Rolle Theorem in Cantor**

In the period from 1874 to 1884, Cantor wrote his main articles on the set theory [38]. In 1872, he began to create his theory of real numbers; in 1874, he proved countability of the set of algebraic numbers; in 1878, he developed the notion of power of set and considered the problem of comparing the power of continuous manifolds of any different dimensions, and came to a paradoxical conclusion that all of them have the same power and are equivalent to a unit segment. "I can see but cannot believe it", wrote Cantor to Dedekind. Cantor concluded that the notion of dimensionality must be based on mutually continuous mapping of varieties on one another.

In 1879, Cantor tried to prove the theorem that two continuous varieties $M_\mu$ and $M_\nu$ with various degrees µ and ν, where $\mu<\nu$, cannot be map one-to-one onto one another continuously. To prove this theorem, Cantor used the theorem of Rolle on the root interval. The remarkable fact is that Cantor referred to Cauchy's *Course of Analysis* of 1821. However, the scheme of his verbal proof was very similar to Bolzano's interpretation of 1817. This attempt to prove the theorem was stated by Cantor in his article entitled "Über einen Satz aus der Theorie der stetigen Mannigfaltigkeiten"(On a theorem from the theory of continuous manifolds) as translated by F.A. Medvedev) [38, p. 36–39]. His proof was not a complete one[20]. However, the theorem on the root interval took a fundamental meaning in his works for the notions of continuity and dimensionality.

## 1886, Weierstrass and substantiation of continuity

In the summer term of 1886 (May-June), Weierstrass lectured on substantiation of the theory of analytical functions. Based on the notion of a limiting point Weierstrass developed the notion of the least upper bound[21] using the theorem on the root interval. Based on this theorem he introduced the notion of connectivity: "proceeding from any point of the continuum, we will always remain there" [36, p. 70].

Many mathematicians of the 19th century addressed these theorems; in many cases, their analysis was pretty interesting.

---

[20] It was L.E.Y. Brauer who provided the first satisfactory proof of the general theorem that the varieties of various number of measurements cannot be mapped one-to-one onto one another continuously and at the same time mutually continuously in his "Math. Ann." of 1910, Bd. 70, s. 161-165.
[21] In doing so Weierstrass used variational methods.

Unfortunately, works of Herman Hankel and many other scientists are beyond the scope of this article.

**Conclusion**

Analysis of the algebraic equation has caused two fundamental statements in the theory of functions (the theorem on the root interval and theorem of the root of a derivative) to appear and the theorem of mean value to be created. In Russian sources, various authors call the theorem on the root interval either 'the second Rolle's theorem' (Shatunovsky [39, p. 121–122]) or 'the theorem of Bolzano–Cauchy' (Shatunovsky's student G.M. Fihtengolz [40, p. 128]). The mean value theorem is called the second theorem of Bolzano–Cauchy [40, p. 131].

It took three hundred years for one of the most fundamental theorems in analysis to form. This theorem is not only of great importance in terms of methodology, but has great applied significance as well, e.g. in differential geometry, functional analysis, mechanics. Luzin proved the theorem on osculating circle and theorem on the center of curvature [41, p. 317] with the help of this theorem. Initially intended for polynomials, Rolle's theorem was extended to continuous functions and enriched their properties. N.N. Luzin said that "this theorem underlies the theoretical development of differential and integral calculus" [41, p. 317].

**References**

1. Больцано Б. Чисто аналитическое доказательство теоремы, что между любыми двумя значениями, дающими результаты противоположного знака, лежит по меньшей мере один действительный корень уравнения / Перевод Э. Кольмана // В кн. Кольман Э. Бернард Больцано. М. – 1955. – С. 170–204. (Bolzano B. Rein analytischer Beweis des Lehrsatzes, dass zwischen je zwey Werthen, die ein entgegengesetzes Resultat gewähren, wenigstens eine reele Wurzel der Gleichung liege. – Prag: Gottlieb Haase. – 1817. – 60 s.)

2. Лопиталь Г. Ф. Анализ бесконечно малых/Пер. с фр. Н.В. Леви под редакцией и со вступительной статьёй А.П. Юшкевича. – М.-Л.: ГТТИ. – 1935 г. – 431 с. (G. de l'Hôpital. Analyse des Infiniment Petits pour l'Intelligence des Lignes Courbes,1696 ). Transl. N. Levi with preface by A.Yushkevich, 1935.

3. Rolle M. Traité d'algèbre, ou principes généraux pour résoudre les questions de mathématique. – Paris. – 1690. – 272 p.

4. Rolle M. Démonstration d'une méthode pour résoudre les égalités de tous les degrés / M. Rolle – Paris. – 1691.

5. Яновская С.А. Мишель Ролль как критик анализа бесконечно малых// Труды института истории естествознания. – М. – 1947 г. – I. – С. 327–346. Yjanovskaja S. Michel Rolle as a critic of infinitesimal analysis // Trudy instituta istorii estestvoznaniia. Moscow, 1947. V.I. – P.327-346.

6. Washington Ch. Michel Rolle and his Method of Cascades // Электронный ресурс: [http://mathdl.maa.org/images/upload_library/46/Washington_Rolle_ed.pdf](http://mathdl.maa.org/images/upload_library/46/Washington_Rolle_ed.pdf)

7. Fontenelle B. Eloge de M. Rolle // Histoire de l`Académie Rouale des Sciences. Paris. – 1719. – P. 94–100.

8. Sinkiewich G. I. Historia dwóch twierdzeń analizy matematycznei: M. Rolle, B. Bolzano, A. Cauchy // Dzieje matematyki polskiei II. Praca zbiorowa pod redakcją Witolda Więsława. Institut Matematyczny Uniwersytetu Wrocławskiego. Wrocław. – 2013. – P. 165 – 181. (Sinkevich G. On the History of two Theorems from Analysis).

9. Мордухай-Болтовской Д. Д. Из прошлого аналитической геометрии. – 1952 г. – Электронный ресурс: [http://wikilivres.ca/wiki/%D0%98%D0%B7_%D0%BF%D1%80%D0%BE%D1%88%D0%BB%D0%BE%D0%B3%D0%BE_%D0%B0%D0%BD%D0%B0%D0%BB%D0%B8%D1%82%D0%B8%D1%87%D0%B5%D1%81%D0%BA%D0%BE%D0%B9_%D0%B3%D0%B5%D0%BE%D0%BC%D0%B5%D1%82%D1%80%D0%B8%D0%B8#cite_ref-10](http://wikilivres.ca/wiki/%D0%98%D0%B7_%D0%BF%D1%80%D0%BE%D1%88%D0%BB%D0%BE%D0%B3%D0%BE_%D0%B0%D0%BD%D0%B0%D0%BB%D0%B8%D1%82%D0%B8%D1%87%D0%B5%D1%81%D0%BA%D0%BE%D0%B9_%D0%B3%D0%B5%D0%BE%D0%BC%D0%B5%D1%82%D1%80%D0%B8%D0%B8#cite_ref-10) (Morduchaj-Boltovskoy. From the past of Analytic Geometry, 1952).

10. Raphson J. Analysis Aequationum universalis. 2 ed. – London. – 1702. – 170 P.

11. Ньютон И. Всеобщая арифметика или книга об арифметическом синтезе и анализе. Пер. А.П. Юшкевича. – М.: Наука. – 1948. – 442 с. (Newton I. Arithmetica Universalis, 1707).

12. Newton I. Method of fluxions. London. – 1786.

13. Simpson T. Essays. London. – 1740. – 144 p.

14. De l'Hospital G.F. Analyse des infiniment petits pour l'intelligence des lignes courbes. – Paris. – 1696.

15. Reynaud Ch.-R. Analyse Démontrée. – Paris. – 1708. – T. I, II. (2-е издание с примечаниями Вариньона. – 1736 г. –Second ed., with Varignon notes).

16. Campbell G. A Method for Determining the Number of Impossible Roots in Adfected Aequations // Philisorhical Transactions of Royal Society. London. – 1727/28. – 35. – P. 515 – 531.

17. A second letter from Mr. Colin Mac Laurin, concerning the Roots of Equations, with the Demonstration of on ther Rulel in Algebra// Phylosophical Transactions of Royal Society. London. – 1729. – 36. – P. 59–96.

18. Clairaut A. C. Élémens d`Algèbre. Paris. – 1746. – 349 p.

19. Эйлер Л. Дифференциальное исчисление. М.-Л.: ГИТТЛ. – 1949. – 580 с. (Euler L. Institutiones calculi differentialis cum eius usu in analysi finitorum ac doctrina serierum, 1755).

20. Lacroix S. F. Élémens d'algèbre, à l'usage de l'Ecole centrale des Quatre-Nations. Paris. – 1830. – 15 ed.

21. Lacroix S.F., Metternich M. Anfangsgründe der Algebra. Mainz – 1811. – 596 s.

22. Kaestner A.G. Angfangsgründe der Analysis endlicher grössen. Göttingen. – 1794. – 590 s.

23. Lagrange J.-L. Traité de la resolution des équations numériques de tous les degrés, avec des notes sur plusieurs points de la théorie des équations algébriques // Oeuvres complètes. Paris. – 1867–1892. – T. 8. – P. 11–367.

24. Rösling Ch. L. Grundlehren von den Formen, Differenzen, Differentialien und Integralien der Functionen. Erlangen. – 1805. – 456 s.

25. Рыхлик К. Теория вещественных чисел в рукописном наследии Больцано // Историко-математические исследования. – 1958 г. – XI. – С. 515–532. (Rychlik K. Theory of real numbers in Bolzano's manuscripts. – Istoriko-matematicheskije issledovanija. Moscow, 1958. – XI. – P.515-532)

26. Cauchy A. Cours d'analyse de l'Ecole royale polytechnique. Première partie: Analyse algébrique // Œuvres complètes. Série 2, tome 3. – Paris. – 1882–1974. – 471 s.

27. Синкевич Г.И. К истории эпсилонтики//Математика в высшем образовании. – 2012. – №10. – С. 149–166.( Sinkevich. On the History of epsilontics//Matematika v vyshem obrasovanii. Nizhny Novgorod, 2012. – 10. – P. 149-166).

28. Drobisch M.W. Grundzüge der Lehre von den höheren Gleichungen. Leipzig. – 1834. – 386 s.

29. Bellavitis G. Sul piu facile modo di trovare le radici reali delle equazioni algebraiche e sopra un nuovo metodo per la determinazione delle radici immaginarie memoria. Venezia. – 1846. – 111 p.

30. Синкевич Г.И. Генрих Эдуард Гейне. Теория функций// Математическое моделирование, численные методы и комплексы программ: межвузовский тематический сборник трудов. Выпуск 18. Под редакцией д-ра физ.- мат. наук, проф. Б.Г. Вагера/ СПбГАСУ. – СПб. – 2012. – С. 6 – 26. (Sinkevich G. Heinrich Eduard Heine. Function Theory//Matematicheskoie modelirovanije, chislennyje metody and complexes of programs. – St.-Petersburg: SPSUACE, 2012. – 18. – P.6-26.

31. Гейне Э. Г. Лекции по теории функций / Э.Г. Гейне. Перевод и примечания Г.И.Синкевич // Математическое моделирование, численные методы и комплексы программ: межвузовский тематический сборник трудов. Выпуск 18. Под редакцией д-ра физ.- мат. наук, проф. Б.Г. Вагера / СПбГАСУ. – СПб. – 2012. – С. 26 – 46. Heine H. Die Elemente der Functionenlehre/Translation and notes Sinkevich// Matematicheskoie modelirovanije, chislennyje metody and complexes of programs. – St.-Petersburg: SPSUACE, 2012. – 18. – P.26-46.

32. Синкевич Г.И. Понятие непрерывности у Дедекинда и Кантора// Труды XI Международных Колмогоровских чтений: сборник статей. – Ярославль: Издательство ЯГПУ – 2013 г. – С. 336–347. (Sinkevich G. A notion of continuity in Dedekind and Cantor's works//Trudy XI Mezhdunarodnych Kolmogorovskich Chtenij. Yaroslavl: YSPU, 2013. – P. 336-347).

33. Karl Weierstrass: Differentialrechnung. Ausarbeitung der Vorlesung an dem König. Gewerbeinstitut zu Berlin im Sommersem / K. Weierstrass, 1861 von H. Schwarz. – in [34].

34. Dugac P. Eléments d'analyse de Karl Weierstrass. Paris. – 1972.

35. Хрестоматия по истории математики. Математический анализ / Под ред. А.П. Юшкевича. – Москва: Просвещение, 1977. – 224 с. (Chrestomatija po istorii matematiki. Matematicheskij Analis. Pod redaktsijej A. Yushkevicha. Moscow, 1977).

36. Weierstrass K. Ausgewählte Kapitel aus der Funktionenlehre. Vorlesung gehalten in Berlin 1886 mit der Akademischen Antrittsrede, Berlin 1857 und drei weiteren Originalarbeiten von K. Weierstrass aus den Jahren 1870 bis 1880/86. Teubner-Archiv für mathematic. Band 9, 272 s. Reprint 1989.

37. Dini U. Fondamenti per la theorica delle funzioni di variabili reali / U. Dini. Pisa. – 1878. – 416 p.

38. Кантор Г. Труды по теории множеств. М.: Наука. – 1985. – 485 с. (Cantor G. Trudy po teorii mnojestv (Works on Set Thery). Moscow, 1985).

39. Шатуновский О.С. Введение в анализ. Одесса. – 1923. – VIII + 244 с. (Shatunovskij O. Vvedenie v analis (Foundation of Analysis). Odessa, 1923.

40. Фихтенгольц Г.М. Основы математического анализа. М. – 1968. – 6 изд. – т. I. – 440 с. Fichtengolts G. Osnovy

matematicheskogo analisa (Foundation of Analysis). Moscow, 1968. - V. I.

41. Лузин Н.Н. Дифференциальное исчисление / Н.Н. Лузин. М.: Советская наука. – 1946. – 451 с. Lusin N. Differntsialnoie ischislenije (Diferential calculus). Moscow, 1946.

42. Stolz, O. B. Bolzano's Bedeutung in der Geschichte der Infinitesimalrechnung / O. Stolz // Mathematische Annalen. – 1881. – Band 18. – P. 255–279.

43. Bolzano B. Purely analytic proof of the theorem that between any two values which give results of opposite sign, there lies at least one real root of the equation; Ewald 1996, pp. 225–248)// Ewald, William B., ed. (1996), From Kant to Hilbert: A Source Book in the Foundations of Mathematics, 2 volumes, Oxford University Press.

44. Robert E. Bradley, C. Edward Sandifer. Cauchy's Course d'analyse. An Annotated Translation. NY: Springer. – 2009. – 432 p.